\documentclass[12pt]{amsart}
\usepackage{graphicx}
\usepackage{amssymb}
\usepackage{amsmath}
\usepackage{amsthm}
\usepackage{amscd}

\newtheorem{thm}{Theorem}[section]
\newtheorem{sub}[thm]{}
\newtheorem{cor}[thm]{Corollary}
\newtheorem{lem}[thm]{Lemma}

\theoremstyle{definition}
\newtheorem{defn}[thm]{Definition}
\theoremstyle{remark}
\newtheorem{rem}[thm]{Remark}
\numberwithin{equation}{section}

\newcommand{\s}{\hfill\blacksquare}

\begin{document}

\title[Monomial Hopf Algebras ]{Monomial Hopf Algebras}
\author[Xiao-wu Chen, Hua-lin Huang, Yu Ye and Pu Zhang
] {Xiao-wu Chen, Hua-lin Huang, Yu Ye, and Pu Zhang$^*$ }

\thanks{$^*$ The corresponding author}
\thanks{Supported in part by the National Natural Science Foundation of China (Grant No. 10271113 and No. 10301033)
and the Europe Commission AsiaLink project  ``Algebras and
Representations in China and Europe$"$ ASI/B7-301/98/679-11}

\maketitle

\vskip10pt

\dedicatory{\begin{center}Dedicated to Claus Michael Ringel \\
on the occasion of his sixtieth birthday\end{center}}
\vskip10pt

\begin{center}
Department of Mathematics \\University of Science and Technology
of China \\Hefei 230026, Anhui, P. R. China
\end{center}

\vskip10pt

\begin{center}
Department of Mathematics\\ Shanghai Jiao Tong University\\
Shanghai 200030, P. R. China
\end{center}
\vskip10pt
\begin{center} xwchen$\symbol{64}$mail.ustc.edu.cn \ \ \ \
hualin$\symbol{64}$ustc.edu.cn\ \ \ \ yeyu$\symbol{64}$ustc.edu.cn

\vskip10pt

pzhang$\symbol{64}$sjtu.edu.cn\end{center} \vskip10pt
\begin{abstract} Let $K$ be a field of characteristic $0$
containing all roots of unity. We classified all the Hopf structures
on monomial $K$-coalgebras, or, in dual version, on monomial
$K$-algebras.
\end{abstract}

\vskip20pt

\centerline {Introduction}

\vskip10pt

In the representation theory of algebras, one uses quivers and
relations to construct algebras, and the resulted algebras are
elementary, see Auslander-Reiten-Smal$\phi$ [ARS] and Ringel
[Rin]. The construction of a path algebra has been dualized by
Chin and Montgomery [CM] to get a path coalgebra. It is then
natural to consider subcoalgebras of a path coalgebra, which are
all pointed.

There are also several works to construct neither commutative nor
cocommutative Hopf algebras via quivers (see e.g. [C], [CR1],
[CR2], [GS]). An advantage for this construction is that a natural
basis consisting of paths is available, and one can relate the
properties of a quiver to the ones of the corresponding Hopf
structures.

\vskip10pt

In [C] Cibils determined all the  graded Hopf structures (with
length grading) on the path algebra $KZ_n^a$ of basic cycle $Z_n$
of length $n$; in [CR1], Cibils and Rosso studied graded Hopf
structures on path algebras; in [GS] E. Green and Solberg studied
Hopf structures on some special quadratic quotients of path
algebras. More recently, Cibils and Rosso [CR2] introduced the
notion of the Hopf quiver of a group with ramification, and then
classified all the  graded Hopf algebras with length grading on
path coalgebras. It turns out that a path coalgebra $KQ^c$ admits
a graded Hopf structure (with length grading) if and only if $Q$
is a Hopf quiver (here a Hopf quiver is not necessarily finite).

\vskip10pt

The cited works above stimulate us to look for finite-dimensional
Hopf algebra structures,  on more quotients of path algebras, or in
dual version, on more subcoalgebras of path coalgebras.

The aim of this paper is to classify  all the Hopf algebra
structures on a monomial algebra, or equivalently, on a monomial
coalgebra.

\vskip10pt

Since a finite-dimensional Hopf algebra is both Frobenius and
co-Frobenius, we first look at the structure of monomial Frobenius
algebras, or dually, the one of monomial coFrobenius coalgebras.
It turns out that each indecomposable coalgebra component of a
non-semisimple monomial coFrobenius coalgebra is $C_d(n)$ with
$d\geq 2$, where $C_d(n)$ is the subcoalgebra of path coalgebra
$KZ_n^c$ with basis the set of paths of length strictly smaller
than $d$. See Section 2.

\vskip10pt

Then by a theorem of Montgomery (Theorem 3.2 in [M2]), a
non-semisimple monomial Hopf algebra $C$ is a crossed product of a
Hopf structure on $C_d(n)$ with a group algebra. Thus, we turn to
study the Hopf structures on $C_d(n)$ with $d\ge 2$. It turns out
that the coalgebra $C_d(n), \ d\geq 2,$ admits a Hopf structure if
and only if $d|n$ (Theorem 3.1). Moreover,
when $q$ runs over primitive $d$-th
roots of unity, the generalized Taft
algebras $A_{n, d}(q)$ gives all the isoclasses of graded Hopf structures on $C_d(n)$
with length grading; while the Hopf structures (not necessarily graded with length grading)
on $C_d(n)$ are exactly the algebras
denoted by $A(n, d, \mu, q)$,  with $q$ a primitive $d$-th root of
unity and $\mu\in K$. These algebras $A(n, d, \mu, q)$ have been studied by Radford [Rad],
Andruskiewitsch and Schneider [AS]. See Theorem 3.6.

\vskip10pt

Note that algebra $A(n, d, \mu, q)$ is given by generators and
relations. In Section 4, we prove that $A(n, d, \mu, q)$ is the
product of $KZ_d^a/J^d$ and $\frac{n}{d}-1$ copies of matrix
algebra $M_d(K)$ when $\mu\ne 0$, and the product of $\frac{n}{d}$
copies of $KZ_d^a/J^d$ when $\mu = 0$, see Theorem 4.3. Hence the
Gabriel quiver and the Auslander-Reiten quiver of $A(n, d, \mu,
q)$ are known.

\vskip10pt

Finally, we introduce the notion of a group datum. By using the
quiver construction of $C_d(n)$, the Hopf structure on it, and
Montgomery's theorem (Theorem 3.2 in [M2]), we get a one to one
correspondence of Galois type between the set of the isoclasses of
non-semisimple monomial Hopf $K$-algebras and the isoclasses of
group data over $K$. This gives a classification of monomial Hopf
algebras.

\section{Preliminaries}

Throughout this paper, $K$ denotes a field of characteristic $0$
containing all roots of unity. By an algebra we mean a
finite-dimensional associative $K$-algebra with identity element.

\vskip10pt

Quivers considered here are always finite. Given a quiver $Q =
(Q_0, Q_1)$ with $Q_0$ the set of vertices and $Q_1$ the set of
arrows, denote by $KQ$, $KQ^a$, and $KQ^c$, the $K$-space with
basis the set of all paths in $Q$, the path algebra of $Q$, and
the path coalgebra of $Q$, respectively. Note that they are all
graded with respect to length grading. For $\alpha\in Q_1$, let $s(\alpha)$ and
$t(\alpha)$ denote respectively the starting and ending vertex of $\alpha$.

\vskip10pt

Recall that the comultiplication of the path coalgebra $KQ^c$ is
defined  by (see [CM])

\begin{align*}
\Delta(p) &=\sum_{\beta\alpha = p}\beta\otimes\alpha \\
&= \alpha_l \cdots \alpha_1\otimes s(\alpha_1)+
\sum_{i=1}^{l-1}\alpha_l \cdots \alpha_{i+1} \otimes
\alpha_{i}\cdots \alpha_1+t(\alpha_l)\otimes\alpha_l \cdots
\alpha_1
\end{align*}
for each path $p = \alpha_l\cdots \alpha_1$ with each $\alpha_i\in
Q_1$; and $\varepsilon(p)=0$ if $l\ge 1$, and $1$ if $l =0$. This
is a pointed coalgebra.

\vskip10pt

Let $C$ be a coalgebra. The set of group-like elements is defined
to be $$G(C): = \{ \ c\in C\ | \ \Delta(c) = c\otimes c, \ \ c
\neq 0\ \}.$$ It is clear $\epsilon(c) = 1$ for $c\in G(C)$. For
$x,y \in G(C)$, denote by
$$P_{x,y}(C): = \{c \in C \ |\ \Delta(c)= c\otimes x + y \otimes
c\},$$ the set of $x,y$-primitive elements in $C$. It is clear
that $\epsilon(c) = 0$ for $c\in P_{x, y}(C)$. Note that
$K(x-y)\subseteq P_{x, y}(C)$. An element $c\in P_{x, y}(C)$ is
non-trivial if $c\notin K(x-y)$.  Note that $G(KQ^c) = Q_0$; and

\vskip10pt

\begin{lem} \ For $x, y\in Q_0$,  we have

$$P_{x, y}(KQ^c) = y(KQ_1)x \oplus K(x-y)$$
where $y(KQ_1)x$ denotes the $K$-space spanned by all arrows from
$x$ to $y$. In particular, there is a non-trivial $x, y$-primitive element
in $KQ^c$ if and only if there is an arrow from $x$ to $y$ in $Q$.
\end{lem}

\vskip10pt

An ideal $I$ of $KQ^a$ is admissible if $J^N\subseteq I \subseteq
J^2$ for some positive integer $N \geq 2$, where $J$ is the ideal
generated by arrows.

An algebra $A$ is elementary if $A/R \cong K^n$ as algebras for
some $n$, where $R$ is the Jacobson radical of $A$. For an
elementary algebra $A$, there is a (unique) quiver $Q$,  and an
admissible ideal $I$ of $KQ^a$, such that $A\cong KQ^a/I$. See
[ARS], [Rin].

\vskip10pt

An algebra $A$ is monomial if there exists an admissible ideal $I$
generated by some paths in $Q$ such that $ A \cong KQ^a/I$.
Dually, we have

\vskip10pt

\begin{defn} \ A subcoalgebra $C$ of $KQ^c$ is called monomial provided
that the following conditions are satisfied

(i) \ \  $C$ contains all vertices and arrows in $Q$;

(ii) \ \ $C$ is contained in subcoalgebra
$C_d(Q):=\bigoplus_{i=0}^{d-1} KQ{(i)}$ for some $d\geq 2$, where
$Q(i)$ is the set of all paths of length $i$ in $Q$;

(iii) \ \ $C$ has a basis consisting of paths.
\end{defn}

\vskip10pt

It is clear by definition that both monomial algebras and monomial
coalgebras are finite-dimensional; and $A$ is a monomial algebra
if and only if the linear dual $A^*$ is a monomial coalgebra.

\vskip10pt

In the following, for convenience, we will frequently pass from a
monomial algebra to a monmial coalgebra by duality. For this we
will use the following

\begin{lem} \ The path algebra $KQ^a$ is exactly the graded dual
of the path coalgebra $KQ^c$, i.e.,
 $$ KQ^a \cong (KQ^c)^{gr};$$
  and for each $d\geq 2$
there is a graded algebra isomorphism :

$$KQ^a/J^d \cong (C_d(Q))^*.$$

\end{lem}

\vskip10pt

\sub{}\rm  \ \ Let $q\in K$ be an $n$-th root of unity. For
non-negative integers $l$ and $m$, the Gaussian binomial
coefficient is defined to be

$${{m+l}\choose l}_q : = \frac{(l+m)!_q}{l!_qm!_q}$$
where

$$l!_q: = 1_q\cdots l_q, \ \ \ \  0!_q: = 1, \ \ \
 l_q: = 1+ q + \cdots q^{l-1}.$$
\vskip10pt

Observe that ${d\choose l}_q = 0$ for $1\le l\le d-1$ if the order
of $q$ is $d$.

\vskip10pt

\sub{}\rm \ \  Denote by $Z_n$ the basic cycle of length $n$, i.e. an oriented
graph with $n$ vertices $e_0, \cdots, e_{n-1}$, and a unique arrow
$\alpha_i$ from $e_i$ to $e_{i+1}$ for each $0\leq i \leq n-1$.
Take the indices modulo $n$.
Denote by $p_i^l$ the path in $Z_n$ of length $l$
starting at $e_i$. Thus we have\  $p_i^0 = e_i$ and \
$p_i^1=\alpha_i$.

\vskip 10pt

For each $n$-th root $q\in K$ of unity, Cibils and Rosso [CR2] have
defined a graded Hopf algebra structure $KZ_n(q)$ (with length
grading) on the path coalgebra $KZ_n^c$ by

$$p_i^l \cdot p_j^m= q^{jl} {{m+l}\choose l} _q p_{i+j}^{l+m},$$
\vskip10pt

\noindent with antipode $S$ mapping $p_i^l$ to $(-1)^l
q^{-\frac{l(l+1)}{2}-il} p_{n-l-i}^l$.

\vskip10pt

 \sub{}\rm \ \ In the following, denote $C_d(Z_n)$ by $C_d(n)$. That is, $C_d(n)$
is the subcoalgebra of $KZ_n^c$ with basis the set of all paths of
length strictly less than $d$.

\vskip10pt

Since ${{m+l}\choose l}_q = 0$ for $m\le d-1, \ l\le d-1, \
l+m\ge d$, it follows that if the order of $q$ is $d$ then
$C_d(n)$ is a subHopfalgebra of $KZ_n(q)$. Denote this graded Hopf
structure on $C_d(n)$ by $C_d(n, q)$.

\vskip10pt

Let $d$ be the order of $q$. Recall that by definition
$A_{n,d}(q)$ is an associative algebra generated by elements $g$
and $x$, with relations

$$ g^n=1, \quad x^d=0, \quad xg=qgx.$$
Then $A_{n,d}(q)$ is a Hopf algebra with comultiplication
$\Delta$, counit $\epsilon$, and antipode $S$ given by

\vskip5pt

\begin{center}
$\Delta(g)=g \otimes g, \quad \quad \quad \epsilon(g)=1,$ \\
\vskip 5pt
$\Delta(x)= x \otimes 1 + g \otimes x, \quad \epsilon(x)=0,$ \\
\vskip 5pt $S(g)=g^{-1}=g^{n-1}, \quad
S(x)=-xg^{-1}=-q^{-1}g^{n-1}x.$
\end{center}

\vskip10pt

In particular, if $q$ is an $n$-th primitive root of unity (i.e.,
$d = n$), then $A_{n,d}(q)$ is the $n^2$-dimensional Hopf algebra
introduced by Taft [T]. For this reason $A_{n,d}(q)$ is called a
generalized Taft algebra in [HCZ].

\vskip10pt

Observe that $C_d(n, q)$ is generated by $e_1$ and $\alpha_0$ as
an algebra. Mapping $g$ to $e_1$ and $x$ to $\alpha_0$, we get a
Hopf algebra isomorphism

$$A_{n,d}(q)\cong C_d(n, q).$$

\vskip 10pt

\sub{}\rm \ \ Let $q\in K$ be an $n$-th root of unity of order
$d$. For each $\mu\in K$, define a Hopf structure $C_d(n, \mu, q)$
on coalgebra $C_d(n)$ by

$$p_i^l\cdot  p_j^m = q^{jl}{m+l \choose l}_q p_{i+j}^{l+m},  \hskip
0.5cm  \mbox{if} \hskip 0.2cm l+m < d, $$ and
$$p_i^l \cdot p_j^m =\mu q^{jl} \frac{(l+m-d)!_q}{l!_qm!_q} ( p_{i+j}^{l+m-d}
-p_{i+j+d}^{l+m-d}), \hskip 0.5cm  \mbox{if} \hskip 0.2cm l+m \geq
d,$$

\vskip10pt

\noindent with antipode

$$S(p_i^l)=(-1)^l q^{-\frac {l(l+1)}{2}-il}p_{n-l-i}^l,$$
\vskip10pt

\noindent where \ $0\leq l, \ m \leq d-1$, and \ $0\leq i, \ j\leq
n-1$. This is indeed a Hopf algebra with identity element $p^0_0 =
e_0$ and of dimension $nd$. Note that this is in general not
graded with respect to the length grading; and that
$$C_d(n, 0, q) = C_d(n, q).$$

\vskip10pt

In [Rad] and [AS] Radford and Andruskiewitsch-Schneider have
considered the following Hopf algebra $A(n, d, \mu, q)$, which as
an associative algebra is generated by two elements $g$ and $x$
with relations

$$g^n=1, \ \ x^d=\mu(1-g^d), \ \ xg=qgx,$$

\vskip10pt

\noindent with comultiplication $\Delta$, counit $\epsilon$, and
the antipode $S$ given as in 1.6.

It is clear that

$$A(n, d, 0, q) = A_{n, d}(q);$$
and if $d=n$ then $A(n, d, \mu, q)$ is the $n^2$-dimensional Taft
algebra.

\vskip 10pt

Observe that $C_d(n, q, \mu)$ is generated by $e_1$ and
$\alpha_0$. By sending $g$ to $e_1$ and $x$ to $\alpha_0$ we
obtain a Hopf algebra isomorphism

$$A(n, d, \mu, q) \cong C_d(n, \mu, q).$$

\vskip10pt

\section{Monomial Forbenius algebras and coFrobenius coalgebras}

\vskip10pt

The aim of this section is to determine the form of monomial
Frobenius, or dually, monomial coFrobenius coalgebras, for later
application. This is well-known, but it seems that there are no
exact references.

\vskip10pt

Let $A$ be a monomial algebra. Thus, $A \cong KQ^a/I$ for a finite
quiver $Q$, where $I$ is an admissible ideal generated by some
paths of lengths $\geq 2$. For $p\in KQ^a$, let $\bar p$ be the
image of $p$ in $A$. Then the finite set
$$\{\ \bar p \in A \ | \ p \ \mbox{does not belong to}\  I\ \}$$
forms a basis of $A$. It is easy to see the following

\vskip10pt

\begin{lem}\ \ Let $A$ be a monomial algebra.  Then

\vskip5pt

(i) \ \ The $K$-dimension of $\operatorname{soc} (Ae_i)$ is the
number of the maximal paths starting at vertex $i$, which do not
belong to $I$.

\vskip5pt

(ii) \ \ The $K$-dimension of $\operatorname{soc} (e_iA)$ is the
number of the maximal paths endding at vertex $i$, which do not
belong to $I$.\end{lem}

\vskip10pt

\begin{lem} \ \ Let $A$ be an indecomposable, monomial algebra.
Then $A$ is Frobenius if and only if $A=k$, or $A \cong
KZ_n^a/J^d$ for some positive integers \ $n$ and $d$, with $d \geq
2$.
\end{lem}

\noindent{\bf Proof}\ \ The sufficiency is straightforward.

If $A$ is Frobenius (i.e. there is an isomorphism $A\cong A^*$ as
left $A$-modules, or equivalently, as right $A$-modules), then the
socle of an indecomposable projective left $A$-module is simple
(see e.g. [DK]). It follows from Lemma 2.1 that there is at most
one arrow starting at each vertex $i$. Replacing "left" by "right"
we observe that there is at most one arrow ending at each vertex
$i$.

On the other hand, the quiver of an indecomposable Frobenius
algebra is a single vertex, or has no sources and sinks (a source
is a vertex at which there is no arrows ending; similarly for a
sink), see e.g. [DK]. It follows that if $A\ne k$ then the quiver
of $A$ is a basic cycle $Z_n$ for some $n$. However it is
well-known that an algebra $KZ_n^a/I$ with $I$ admissible is
Frobenius if and only if $I=J^d$ for some $d \geq 2$. $\s$

\vskip10pt

The dual version of Lemma 2.2 gives
the following

\vskip10pt

\begin{lem} \ \ Let $A$ be an indecomposable, monomial coalgebra.
Then $A$ is coFrobenius (i.e., $A^*$ is Frobenius) if and only if
$A=k$, or $A \cong C_d(n)$ for some positive integers \ $n$ and
$d$, with $d \geq 2$.
\end{lem}

\vskip10pt

An algebra $A$ is called Nakayama, if each indecomposable
projective left and right module has a unique composition series.
It is well-known that an indecomposable elementary algebra is
Nakayama if and only if its quiver is a basic cycle or a linear
quiver $A_n$ (see [DK]). Note that a finite-dimensional Hopf
algebra is Frobenius and coFrobenius (see for e.g. [M1], p.18).

\vskip10pt

\begin{cor} \ An algebra is a monomial Frobenius algebra
if and only if it is elementary Nakayama Frobenius. Hence, a Hopf
algebra is monomial if and only if it is  elementary and
Nakayama.
\end{cor}

\vskip10pt

\section{Hopf structures on coalgebra $C_d(n)$}

\vskip 10pt The aim of this section is to give a numerical
description such that coalgebra $C_d(n)$ admits Hopf structures
(Theorem 3.1), and then classify all the  (graded, or not
necessarily graded) Hopf structures on $C_d(n)$ (Theorem 3.6).

\vskip10pt

\begin{thm}\ Let $K$ be a field of characteristic $0$, containing an $n$-th primitive root of unity.
Let $d\ge 2$ be a positive integer. Then coalgebra $C_d(n)$ admits
a Hopf algebra structure if and only if $d|n$.
\end{thm}

The sufficiency follows from 1.6, or 1.7. In order to prove the
necessity we need some preparations.

\vskip10pt

\begin{lem}\ Suppose that the coalgebra $C_d(n)$ admits a Hopf algebra structure.
Then

\vskip5pt

(i) \ \ The set $\{e_0, \cdots, e_{n-1}\}$ of the vertices in
$C_d(n)$ forms a cyclic group, say, with identity element $1=e_0$.
Then $e_1$ is a generator of the group.

\vskip5pt

(ii) \ \ Set $g: = e_1$. Then up to a Hopf algebra isomorphism we
have for any  $i$ such that $0\le i\le n-1$

$$\alpha_{i} \cdot g = q\alpha_{i+1} + \kappa_{i+1}(g^{i+1} - g^{i+2})$$
and

$$g\cdot \alpha_{i}= \alpha_{i+1} + \lambda_{i+1}(g^{i+1} - g^{i+2}),$$
where $q, \lambda_i, \kappa_i\in K$, with $q^n = 1$.
\end{lem}

\noindent {\bf Proof} \ \ (i) \ \ Since $C_d(n)$ is a Hopf
algebra, it follows that $G(C_d(n)) = \{e_0, \cdots, e_{n-1}\}$ is
a group, say with identity element $e_0$.  Since $\alpha_0$ is a
non-trivial $e_0, e_1$-primitive element, it follows that
$\alpha_0e_1$ is a non-trivial $e_1, e_1^2$-primitive element,
i.e., there is an arrow in $C_d(n)$ from $e_1$ to $e_1^2$. Thus
$e_1^2 = e_2$. A similar argument shows that $e_i = e_1^i$ for any
$i$.

\vskip10pt

(ii) \ \  Since both $\alpha_i g$ and $g\alpha_i$ are non-trivial
$g^{i+1}, g^{i+2}$-primitive elements, it follows that

$$\alpha_i\cdot g = w_{i+1} \alpha_{i+1} + \kappa_{i+1} (g^{i+1} -g^{i+2})$$
and

$$g\cdot \alpha_i = y_{i+1} \alpha_{i+1} + \lambda'_{i+1} (g^{i+1} -g^{i+2})$$
with $w_i, \kappa_i, y_i, \lambda'_i \in K.$

\vskip 5pt

Since $g^n\cdot \alpha_0 = \alpha_0$, it follows that $y_1\cdots
y_n = 1$.  Set $\theta_j: = y_{j+1}\cdots y_n, \ 1\le j \le n-1$,
and $\theta_n: = 1$. Define a linear isomorphism $\Theta: \ C_d(n)
\longrightarrow C_d(n)$ by

\begin{align*} p_i^l&\longmapsto
(\theta_i\cdots \theta_{i+l-1}) p_i^l.
\end{align*}
In particular $\Theta(e_i) = e_i$ and $\Theta(\alpha_i) =
\theta_i\alpha_i$. Then $\Theta: \ C_d(n) \longrightarrow C_d(n)$
is a coalgebra map.  Endow $C_d(n) = \Theta(C_d(n))$ with the Hopf
algebra structure via the given Hopf algebra structure of $C_d(n)$
and $\Theta$. Then in $\Theta(C_d(n))$ we have

\begin{align*}
g\cdot(\theta_i\alpha_i)&= \Theta(g)\cdot\Theta(\alpha_i) =
\Theta(g\cdot\alpha_i) \\& = y_{i+1}\Theta(\alpha_{i+1})
+\lambda'_{i+1}(g^{i+1}-g^{i+2}) \\&=
y_{i+1}\theta_{i+1}\alpha_{i+1}
+\lambda'_{i+1}(g^{i+1}-g^{i+2}).\end{align*} Since $\theta_i =
y_{i+1}\theta_{i+1}$, it follows that in $\Theta(C_d)$ we have

$$g\cdot\alpha_i = \alpha_{i+1} + \lambda_{i+1}(g^{i+1} - g^{i+2})$$
(with $\lambda_{i+1} = \frac{\lambda'_{i+1}}{\theta_i}$). Assume
that now in $\Theta(C_d(n))$ we have
$$\alpha_i \cdot g = q_{i+1} \alpha_{i+1} + \kappa_{i+1} (g^{i+1} -g^{i+2})$$
Since $\alpha_0 g^n=\alpha_0$, it follows that $q_1\cdots q_n=1$.
However, $(g\cdot\alpha_i)\cdot g=g\cdot(\alpha_i\cdot g)$ implies
$q_i=q_{i+1}$ for each $i$. Write $q_i=q$. Then $q^n=1$. This
completes the proof. $\blacksquare$

\vskip 10pt

\begin{lem} \ \ Suppose that there is a Hopf algebra structure on $C_d(n)$. Then
up to a Hopf algebra isomorphism we have

$$p_i^l \cdot p_j^m \equiv q^{jl}{m+l \choose l}_q\  p_{i+j}^{l+m}\hskip0.5cm (mod \ \ C_{l+m}(n))$$
for $0\leq i, j \leq n-1$, and for \ $l, m\ \le d-1$, where $q\in
K$ is an $n$-th root of unity.
\end{lem}

\vskip 5pt

\noindent {\bf Proof} \ Use induction on $N: = l+m$. For $N=0$ or
$1$, the formula follows from  Lemma 3.2. Assume that the formula
holds  for $N\leq N_0-1$. Then for $N=N_0\ge 1$ we have

\begin{align*}
\bigtriangleup(p_i^l \cdot p_j^m)
&=\bigtriangleup(p_i^l) \cdot  \bigtriangleup(p_j^m) \\
&=(\sum_{r=0}^{l}p_{i+r}^{l-r}\otimes p_i^r )\cdot
(\sum_{s=0}^{m}p_{j+s}^{m-s}\otimes p_j^s ) \\
&=\sum_{k=0}^{N_0} \ \sum_{r+s =k, 0\le r\le l, 0\le s\le
m}p_{i+r}^{l-r}p_{j+s}^{m-s}\otimes p_i^rp_j^s
\\& = p_i^lp_j^m \otimes g^{i+j}+ g^{i+j+N_0}\otimes p_i^lp_j^m
\\&+ \sum_{k=1}^{N_0-1} \ \sum_{r+s =k, 0\le r\le l, 0\le s\le
m}p_{i+r}^{l-r}p_{j+s}^{m-s}\otimes p_i^rp_j^s.
\end{align*}
By the induction hypothesis  for each $r$ and $s$ with $1\le k: =
r+s\le N_0-1$ we have

$$p_i^r \cdot p_j^s \equiv q^{jr}{k \choose r}_q\  p_{i+j}^{k}
\hskip0.5cm (mod \ C_{k}(n))$$ and

$$p_{i+r}^{l-r}\cdot  p_{j+s}^{m-s} \equiv q^{(j+s)(l-r)}{N_0-k \choose l-r}_q\
p_{i+j+k}^{N_0-k}\hskip0.5cm (mod  \ C_{N_0-k}(n)).$$ It follows
that

$$\bigtriangleup(p_i^l \cdot p_j^m) \equiv p_i^l\cdot p_j^m\otimes g^{i+j} \  + \
 g^{i+j+N_0}\otimes p_i^l\cdot p_j^m + \Sigma\hskip0.4cm (mod \
\bigoplus_{1\le k\le N_0-1} C_{N_0-k}(n)\otimes C_{k}(n))$$ where

\begin{align*}\Sigma &= q^{jl}\sum_{k=1}^{N_0-1}
\ \sum_{r+s =k,0\le r\le l, 0\le s\le m} q^{sl-sr}{k\choose
r}_q{N_0-k \choose l-r}_q p_{i+j+k}^{N_0-k}\otimes p_{i+j}^k
\\& = q^{jl}\sum_{k=1}^{N_0-1} {N_0 \choose l}_q p_{i+j+k}^{N_0-k}\otimes
p_{i+j}^k.
\end{align*}
Note that  in the last equality  the  following identity has been
used (see e.g. Proposition IV.2.3 in [K]):

$$\sum_{r+s=k} q^{sl-sr}{k
\choose r}_q {N_0-k \choose l-r}_q ={N_0 \choose l}_q, \ \
0<k<N_0.$$

\vskip10pt

Now, put $X:=p_i^l  p_j^m - q^{jl}{N_0 \choose l}_q
p_{i+j}^{N_0}$. Then by the computation above we have \vskip10pt

$$\bigtriangleup(X) \equiv X \otimes g^{i+j}+  g^{i+j+N_0} \otimes X
\hskip 0.5cm (mod \bigoplus_{1\le k \leq N_0-1}
C_{N_0-k}(n)\otimes C_k(n)). $$

\vskip10pt

\noindent Let $X= \sum_{v\ge 0} c_v$, where $c_v$ is the $v$-th
homogeneous component with respect to the length grading. Then we
have

$$\sum_v\bigtriangleup(c_v) \equiv \sum_{v}( c_v\otimes g^{i+j}+ g^{i+j+N_0}\otimes
c_v) \hskip 0.5cm (mod \ \bigoplus_{1\le k \leq N_0-1}
C_{N_0-k}(n)\otimes C_{k}(n)). $$ Since the elements in
$C_{N_0-k}(n)\otimes C_{k}(n)$ are
 of degrees strictly smaller than $N_0$, it follows that for $v\ge
 N_0$ we have

\vskip10pt

$$\bigtriangleup(c_v) = c_v\otimes g^{i+j}+ g^{i+j+N_0}\otimes
c_v.$$

\vskip10pt

\noindent Now for each $v\ge N_0\ge 1$, \ note that in the right
hand side of the above equality terms are of degree $(v, 0)$ or
$(0, v)$; but in the left hand side if $c_v\ne 0$ then it really
contains a term of degree which is neither $(v, 0)$ nor $(0, v)$.
This forces $c_v = 0$ for $v\ge N_0$. It follows that

$$p_i^l  p_j^m = q^{jl}{N_0 \choose l}_q p_{i+j}^{N_0} + X
\equiv q^{jl}{N_0 \choose l}_q p_{i+j}^{N_0}\hskip 0.5cm (mod \ \
C_{N_0}(n)).$$ This completes the proof. $\blacksquare$

\vskip 10pt

By a direct analysis from the definition of the Gaussian binomial
coefficients we have

\vskip5pt

\begin{lem} \ Let $1\ne q\in K$ be an $n$-th root of unity of
order $d$. Then $${m+l \choose l}_q = 0 \ \ \text{if and only if}
\ \ \big[\frac{m+l}{d} \big]-\big[ \frac{m}{d}\big]-\big[
\frac{l}{d} \big]>0,$$ where $[x]$ means the integer part of $x$.
\end{lem}

\vskip10pt

\sub{}\rm {\bf Proof of Theorem 3.1:} \ \ Assume that $C_d(n)$
admits a Hopf algebra structure. Let $q$ be the $n$-th root of
unity as appeared in Lemma 3.3 with order $d_0$. It suffices to
prove $d=d_0$. Since $C_d(n)$ has a basis $p_i^l$ with $l\le d-1$
and $0\le i\le n-1$,  it follows from Lemma 3.3 that

$${m+l \choose l}_q=0 \ \ \text{for} \ \ l,\  m\leq d-1, \ \ l+m \ge d.$$
While by Lemma 3.4
$${m+l \choose l}_q = 0 \ \text{if and only if} \ \big[\frac{m+l}{d_0}
\big]-\big[ \frac{m}{d_0}\big]-\big[ \frac{l}{d_0} \big]>0.$$

\vskip10pt

(Note that here we have used the assumption that $K$ is of
characteristic $0$: since $K$ is of characteristic zero, it
follows that ${m+l \choose l}_1$ can never be zero. Thus $q\ne 1$,
and then Lemma 3.4 can be applied.)

\vskip10pt

Take $l=1$ and $m=d-1$. Then we have
$[\frac{d}{d_0}]-[\frac{d-1}{d_0}]>0$. This means $d_0|d$. Let
$d=kd_0$ with $k $ a positive integer. If $k>1$, then by taking
$l=d_0 $ and $m=(k-1)d_0$ we get a desired contradiction ${l+m
\choose l}_q \neq 0$. $\s$

\vskip10pt

\begin{thm}\ Assume that $K$ is a field of characteristic $0$, containing an $n$-th primitive root of unity.
Let $d|n$ with $d\ge 2$. Then

\vskip5pt

(i)\ \ Any graded Hopf structure (with length grading) on $C_d(n)$
is isomorphic to (as a Hopf algebra) some $C_d(n, q)\cong A_{n,
d}(q)$, where $C_d(n, q)$ and $A_{n, d}(q)$ are given as in 1.6.

\vskip5pt

(ii) \ \ Any Hopf strucutre (not necessarily graded) on $C_d(n)$ is isomorphic to (as a
Hopf algebra) some $C_d(n, \mu, q)\cong A(n, d, \mu, q)$, where
$C_d(n, \mu, q)$ and $A(n, d, \mu, q)$ are given as in 1.7.

\vskip5pt

(iii)\ \ If $A(n_1, d_1, \mu_1, q_1) \simeq A(n_2, d_2, \mu_2,
q_2) $ as Hopf algebras, then $n_1 = n_2, d_1 = d_2$, $q_1=q_2$.

If $d \ne n$, then $A(n, d, \mu_1, q) \simeq A(n, d, \mu_2, q) $
as Hopf algebras if and only if $\mu_1=\delta^d \mu_2$ for some
$0\ne \delta \in K$, and $A(n, n, \mu_1, q) \simeq A(n, n, \mu_2,
q)$ for any $\mu_1, \mu_2\in K$.

\vskip5pt

In particular, for each $n$, $C_{d}(n, q_1)$ is isomorphic to
$C_{d}(n, q_2)$ if and only if $q_1 = q_2$.
\end{thm}

\vskip5pt

\noindent {\bf Proof} \quad (i) \ \ By Lemma 3.3 and by the proof
of Theorem 3.1 we see that any graded Hopf algebra on $C_d(n)$ is
isomorphic to $C_d(n, q)$ for some root $q$ of unity of order $d$.

\vskip5pt

(ii) \ Assume that $C_d(n)$ is a Hopf algebra. By Lemma 3.2 we
have

$$\alpha_0\cdot e_1 = q e_1\cdot \alpha_0 + \kappa (e_1-e_1^2)$$
for some primitive $d$-th root $q$. Set $X: = \alpha_0 +
\frac{\kappa}{q-1} (1-e_1)$. Then $Xe_1=qe_1X$. Since
$\Delta(X)=e_1 \otimes X + X \otimes 1$, it follows that

\begin{align*}
\Delta(X^d)&=(\Delta(X))^d \\&
= \sum_{i=0}^{d} {d \choose i}_q e_{d-i} X^i \otimes X^{d-i}\\
&= e_d\otimes X^d + X^d\otimes 1,
\end{align*}
where in the last equality we have used the fact that ${d \choose
i}_q = 0$ for $1\le i\le d-1$. Since there is no non-trivial $1,
e_d$-primitive element in $C_d(n)$, it follows that  $X^d=\mu
(1-e_1^d)$ for some $\mu \in K$. Hence we obtain an algebra map
\begin{center}
$F: A(n, d, \mu, q) \longrightarrow C_d(n)$
\end{center}
such that $F(g)= e_1$ and $F(x)=X$. Since $C_d(n)$ is generated by
$e_1$ and $\alpha_0$ by Lemma 3.3, it follows that $F$ is
surjective, and hence an algebra isomorphism by comparing the
$K$-dimensions. It is clear that $F$ is also a coalgebra map,
hence a bialgebra isomorphism, which is certainly a Hopf
isomorphism ([S]).

\vskip5pt

(iii)\ \ If $C_{d_1}(n_1, \mu_1, q_1)\cong C_{d_2}(n_2, \mu_2,
q_2)$, then their groups of the group-like elements are
isomorphic. Thus $n_1 = n_2$, and hence $d_1 = d_2$ by comparing
the $K$-dimensions. The remaining assertions can be  easily
deduced. We omit the details.  $\blacksquare$

\vskip10pt

\begin{rem}\ The following example shows that, the assumption "$K$ is of characteristic $0$"
is really needed in Theorem 3.1.

\vskip10pt

Let K be a field of characteristic 2, and let $n\ge 2$ be an arbitrary
integer. Then each graded Hopf algebra structure on $C_2(n)$ is
given by (up to a Hopf algebra isomorphism):

\begin{center} $g^j
\alpha_i=\alpha_i g^j=\alpha_{i+j}, \ \ \ \alpha_i \alpha_j=0,$
\end{center}

\begin{center}
$S(\alpha_i)=\alpha_{n-i-1}, \ \ \ S(g^j)=g^{n-j}$
 \end{center}
 for all $0\leq i, j \leq n-1$.

\vskip10pt

(In fact, consider the Hopf algebra structure $KZ_n(1)$ on $Z_n$.
Its subcoalgebra $C_2(n)$ is also a subalgebra, which is exactly
the given Hopf algebra. On the other hand, for each graded Hopf
algebra over $C_2(n)$, the corresponding $q$ in Lemma 3.3 must
satisfy ${2 \choose 1}_q= 1+q = 0$,  and hence $q=1$. Then the
assertion follows from Lemma 3.3.) \end{rem}

\begin{rem}\ \ It is easy to determine the automorphism group
of the Hopf algebra $A(n, d, \mu, q)$: it is $K-\{0\}$ if $\mu =
0$ or $d = n$, and $Z_d$ otherwise.
\end{rem}

\vskip 10pt

\section{The Gabriel quiver and the Auslander-Reiten quiver
of $A(n, d, \mu, q)$}

\vskip10pt

The aim of this section is to determine the Gabriel quiver and the
Auslander-Reiten quiver of algebra $A(n, d, \mu, q)\cong C_d(n,
\mu, q)$, where $q$ is an $n$-th root of unity of order $d$.

\vskip 10pt

We start from the central idempotent decomposition of  $A: =A(n,
d, \mu, q)$.

\vskip10pt

\begin{lem} \ The center of $A$ has a linear basis $\{1, g^d,
g^{2d}, \cdots, g^{n-d}\}$.

Let $\omega\in K$ be a root of unity of order $\frac{n}{d}$.  Then we have
the central idempotent decomposition \ $1=c_0+c_1+ \cdots +c_t$ \
with $c_i= \frac{d}{n}\sum_{j=0}^t (\omega^i g^d)^j$ for all $0
\leq i \leq t$, where $t=\frac{n}{d}-1$.
\end{lem}

\vskip5pt

\noindent {\bf Proof}\ \ By 1.7 the dimension of $A$ is $nd$, thus
$\{g^ix^j|0 \leq i \leq n-1,0\leq j \leq d-1\}$ is a basis
of $A$. An element $c=\sum a_{ij} g^i x^j $ is in the center of
$A$ if and only if $xc=cx$ and $gc=cg$. By comparing the
coefficients, we get $a_{ij}=0$ unless $j=0$ and $d|i$. Obviously,
$g^d$ is in the center. It follows that the center of $A$ has a
basis $\{1, g^d, g^{2d}, \cdots, g^{n-d}\}$.

Since $\sum_{i=0}^t (\omega^j)^i = 0$ for each $1\le j\le t$, it
follows that

\begin{align*} &c_0+c_1+ \cdots +c_t=\frac{d}{n}\sum_{j=0}^tg^{dj}\sum_{i=0}^t
(\omega^j)^i =\frac{d}{n}(\sum_{i=0}^t 1 + \sum_{j=1}^tg^{dj}) \\&
=\frac{d}{n}(t+1) = 1;
\end{align*}
and
\begin{align*} c_ic_i'
&=\frac{d^2}{n^2}\sum_{0\le j, j'\le t}
g^{d(j+j')}\omega^{ij+i'j'}
\\& = \frac{d^2}{n^2}\sum_{k=0}^{2t}
g^{dk}\omega^{i'k} \sum_{0\le j\le \mbox{min}\{k,t\}, 0\le k-j\le
t}\omega^{(i-i')j}
\\& =\frac{d^2}{n^2}
(\sum_{k=0}^{t}g^{dk}\omega^{i'k} \sum_{0\le j\le
k}\omega^{(i-i')j}+\sum_{k=t+1}^{2t}
g^{dk}\omega^{i'k}\sum_{k-t\le j\le t}\omega^{(i-i')j})
\\& =\frac{d^2}{n^2}(\sum_{k=0}^{t}
g^{dk}\omega^{i'k} \sum_{0\le j\le
k}\omega^{(i-i')j}+\sum_{k'=0}^{t-1}
g^{dk'}\omega^{i'k'}\sum_{1+k'\le j\le t}\omega^{(i-i')j})
\\& =\frac{d^2}{n^2}(g^{dt}\omega^{i't}\sum_{0\le j\le
t}\omega^{(i-i')j}+\sum_{k=0}^{t-1} g^{dk}\omega^{i'k}\sum_{0\le
j\le t}\omega^{(i-i')j})\\&
=\frac{d^2}{n^2}(g^{dt}\omega^{i't}\delta_{i,i'}(t+1)+\sum_{k=0}^{t-1}
g^{dk}\omega^{i'k}\delta_{i,i'}(t+1))\\& =
(t+1)\frac{d^2}{n^2}\delta_{i,i'}\sum_{k=0}^{t}
g^{dk}\omega^{i'k}\\& = \delta_{i,i'}c_i
\end{align*}
where $\delta_{i, i'}$ is the Kronecker symbol. This completes the
proof. $\blacksquare$

\vskip 10pt

\begin{lem}
Let $B=B(d,\lambda,q)$ be an algebra generated by $g$ and $x$ with
relations $\{g^d=1, x^d=\lambda, xg=qgx \}$, where $\lambda, q \in
K$, and $q$ is a root of unity of order $d$.

\vskip5pt

(i)\ \  If $\lambda =0$, then $B\simeq KZ_n^d/{J^d}$.

\vskip5pt

(ii) \ \ If $\lambda \neq 0$, then $B\simeq M_d(K)$.
\end{lem}

\vskip 5pt

\noindent {\bf Proof}\ \ (i)\ \ Note that if $\lambda=0$, then $B
\simeq A(d, d, 0, q)\cong C_d(d, 0, q)$, which is a
$d^2$-dimensional Taft algebra. By the self-duality of the Taft
algebras (see [C], Proposition 3.8) we have algebra isomorphisms
 $$B \cong A(d,d,0,q)
\simeq A(d, d, 0, q)^* \simeq C_d(d, 0, q)^* \simeq
KZ_d^a/{J^d}.$$

\vskip5pt

(ii) \ \ If $\lambda \neq 0$, then define an algebra homomorphism
$\phi: B \longrightarrow M_d(K)$:

 $$\phi(g)=
\begin{pmatrix}
  1 &  &  &  &  \\
    & q &  &  &  \\
    &  & q^2 &  &  \\
    &  & & \ddots &  \\
    &  &  &  & q^{d-1}
\end{pmatrix}$$
and

$$\phi(x)=
\begin{pmatrix}
  0 &  1  &  &  \\
  &  0 & 1   &  &  \\
    & & \ddots &  &  \\
    &  & \ddots &  & &0 &1   \\
    \lambda & &  & & & & 0
\end{pmatrix}.$$
Note that $\phi$ is well-defined. It is easy to check that
$\phi(g)$ and $\phi(x)$ generate the algebra $M_d(K)$. Thus $\phi$
is a surjective map. However, the dimension of $B$ is at most
$d^2$, thus $\phi$ is an algebra isomorphism. $\s$

\vskip 10pt

Now we are ready to prove the main result of this section.

\vskip 10pt

\begin{thm}\ \ Write $A=A(n,d,\mu,q)$ and $t=\frac{n}{d}-1$.

\vskip5pt

(i)\ \  If $\mu \neq 0$, then $A \simeq KZ_d^a/{J^d} \times M_d(K)
\times \cdots \times M_d(K)$ (with $t$ copies of
$M_d(K)$).

\vskip5pt

(ii)\ \  If $\mu=0$, then $A \simeq KZ_d^a/{J^d} \times
KZ_d^a/{J^d}\times \cdots \times KZ_d^a/{J^d}$ (with $\frac{n}{d}$ copies
of $KZ_d^a/{J^d}$).
\end{thm}

\vskip 5pt

\noindent {\bf Proof}\ \ By Lemma 4.1 we have $A\cong c_0A \times
c_1A\times \cdots \times c_tA$ as algebras. Write $A_i=c_iA$. Note
that $c_ig^d=\omega ^{-i} c_i$ for all $0\leq i \leq t$. It
follows that $\{c_i g^kx^j \ | \ 0 \leq k \leq d-1, \ 0\leq j \leq
d-1 \}$ is a linear basis of $A_i$. Let $\omega_0\in K$ be an
$n$-th primitive root of unity such that $\omega_0^d = \omega$.
Obviously, as an algebra each $A_i$ is generated by $\omega_0^i
c_ig$ and $c_ix$, satisfying

$$(\omega_0^ic_ig)^d=c_i, \hskip 0.6cm (c_ix)^d=c_i
\mu(1-g^d)=c_i\mu (1-\omega^{-i})$$ and
$$(c_ix)(\omega_0^ic_ig)=q(\omega_0^ic_ig)(c_ix).$$
Note that $c_i$ is the identity of $A_i$. Thus we have an algebra
homomorphism
$$\theta_i: \ B(d, \mu(1-\omega^{-i}), q)\longrightarrow A_i$$
such that $\theta_i(g)=\omega_0^i c_ig$ and $\theta_i(x)=c_ix$. A
simple dimension argument shows that $\theta_i$ is an algebra
isomorphism. Note that $\mu (1-\omega^{-i})=0$ if and only if
$\mu=0$ or $i=0$. Then the assertion follows from Lemma 4.2.
$\blacksquare$

\vskip 10pt

\begin{cor} \ \ The Gabriel
quiver of algebra $A(n, d, \mu, q)$ is the disjoint union of  a
basic d-cycle and $t$ isolated vertices if $\mu\ne 0$, and  the
disjoint union of $\frac{n}{d}$ basic d-cycles if $\mu = 0$. \end{cor}

\vskip 10pt

Since the Auslander-Reiten quiver of $\Gamma (KZ_d^a/{J^d})$ is
well-known (see e.g. [ARS], p.111), it follows that the
Auslander-Reiten quiver of $A(n, d, \mu, q)$ is clear.

\vskip10pt

\section{Hopf structures on monomial algebras and coalgebras}

\vskip10pt

The aim of is section is to classify non-semisimple monomial Hopf $K$-algebras, by
establishing an one to one correspondence
between the set of
the isoclasses of non-semisimple monomial Hopf
$K$-algebras and
the isoclasses of group data over $K$.

\vskip10pt

\begin{thm} (i) \ \ Let $A$ be a monomial algebra. Then $A$ admits a Hopf
algebra structure if and only if $A\cong k\times \cdots\times k$
as an algebra, or
$$A\cong KZ_n^a/J^d \times \cdots \times KZ_n^a/J^d
$$ as an algebra, for some $d \geq 2$ dividing $n$.

\vskip5pt

(ii) \ \ Let $C$ be a monomial coalgebra. Then $C$ admits a Hopf
algebra structure if and only if $C\cong k\oplus\cdots\oplus k$ as
a coalgebra, or
$$C\cong C_d(n) \oplus\cdots \oplus C_d(n)$$ as a
coalgebra, for some $d \geq 2$ dividing $n$.
\end{thm}

\vskip5pt

{\bf Proof} \ \ By duality it suffices to prove one of them. We
prove (ii).

\vskip5pt

If $C = C_1 \oplus\cdots \oplus C_l$ as a coalgebra, where each
$C_i\cong C_1$ as coalgebras, and $C_1$ admits Hopf structure
$H_1$, then $H_1\otimes KG$ is a Hopf structure on $C$, where $G$
is any group of order $l$. This gives the sufficiency.

\vskip10pt

Let $C$ be a monomial coalgebra admitting a Hopf structure. Since
a finite-dimensional Hopf algebra is coFrobenius, it follows from
Lemma 2.3 that as a coalgebra $C$ has the form $C= C_1 \oplus
\cdots \oplus C_l$ with each $C_i$ indecomposable as coalgebra,
and $C_i = k$ or $C_i = C_{d_i}(n_i)$ for some $n_i$ and \
$d_i\geq 2$.

\vskip10pt

We claim that if there exists a $C_i = k$, then $C_j =k$ for all
$j$. Thus, if $C\ne k\oplus \cdots \oplus k$, then $C$ is of the
form
$$C = C_{d_1}(n_1) \oplus
\cdots \oplus C_{d_l}(n_l)$$ as a coalgebra, with each $d_i\geq
2$.

(Otherwise, let $C_j = C_{d}(n)$ for some $j$. Let $\alpha$ be an
arrow in $C_j$ from $x$ to $y$. Let $h$ be the unique group-like
element in $C_i = k$. Since the set $G(C)$ of the group-like
elements of $C$ forms a group, it follows that there exists an
element $k\in G(C)$ such that $h = kx$. Then $k\alpha$ is a $h,
ky$-primitive element in $C$. But according to the coalgebra
decomposition $C = C_1 \oplus \cdots \oplus C_l$ with $C_i = Kh$,
$C$ has no $h, ky$-primitive elements. A contradiction.)

\vskip10pt

Assume that the identity element $1$ of $G(C)$ is contained in
$C_1 = C_{d_1}(n)$. It follows from a theorem of Montgomery (
Theorem 3.2, [M2]) that $C_1$ is a subHopfalgebra of $C$, and that
$$g_i^{-1}C_{d_i}(n_i) = C_{d_i}(n_i)g_i^{-1}= C_{d_1}(n_1)$$ for
any $g_i\in G(C_{d_i}(n_i))$ and for each $i$. By comparing the
numbers of group-like elements in $g_i^{-1}C_{d_i}(n_i)$ and in
$C_{d_1}(n_1)$ we have $n_i = n_1 = n$ for each $i$.  While by
comparing the $K$-dimensions we see that $d_i = d_1 = d$ for each
$i$. Now, since $C_1 = C_{d}(n)$ is a Hopf algebra, it follows
from Theorem 3.1 that $d$ divides $n$. $\s$

\vskip10pt

\sub{}\rm \ For convenience, we call a Hopf structure on a
monomial coalgebra $C$ a monomial Hopf algebra. Note that a monomial Hopf algebra
is not necessarily graded with length grading, by Lemma (iii) below.

\vskip10pt

\noindent {\bf Lemma.} \  \ {\it Let $C$ be a non-semisimple, monomial
Hopf algebra.

\vskip5pt

(i) \ \ Let $C_1$ be the indecomposable coalgebra component
containing the identity element $1$.  Then $G(C_1)$ is a cyclic
group contained in the center of $G(C)$.

\vskip5pt

(ii) \ \ There exists a unique element $g\in C$ such that there is
a non-trivial $1, g$-primitive element in $C$. The element $g$ is
a generator of $G(C_1)$.

\vskip5pt

(iii) \ \ As an algebra, $C$ is generated by $G(C)$ and a
non-trivial $1, g$-primitive element $x$, satisfying

$$x^d=\mu (g^d-1)$$

\vskip5pt

\noindent for some $\mu \in K$, where $d
=\frac{\operatorname{dim}_KC_1}{o(g)}$, $o(g)$ is the order of
$g$.

\vskip5pt

(iv) \ \  There exists a one-dimensional $K$-representation $\chi$
of $G$ such that

$$x\cdot h = \chi(h)h\cdot x, \ \ \forall h \in G,$$
and  $\mu = 0$ if $o(g) = d$ (note that $d = o(\chi(g))$); and
$\chi^d=1$ if $\mu \neq 0$ and $g^d \neq 1$.}

\vskip5pt

\noindent {\bf Proof} \ \ (i) \ \  Note that $C_1$ is a
subHopfalgebra of $C$ by Theorem 3.2 in [M2]. By Theorem 5.1(ii)
we have $C_1\cong C_d(n)$ as a coalgebra. It follows from Lemma
3.3 that $G(C_1)$ is a cyclic group. By Theorem 5.1(ii) we can
identify each indecomposable coalgebra component $C_i$ of $C$ with
$C_d(n)$. For any $h\in G(C)$ with $h\in C_i$, note that
$h\alpha_0$ is a non-trivial $h, he_1$-primitive element in $C_i$,
and $\alpha_0h$ is a non-trivial $h, e_1h$-primitive element in
$C_i$. This implies that there is an arrow in $C_i = C_{d}(n)$
from $h$ to $he_1$, and that there is an arrow in $C_i$ from $h$
to $e_1h$. Thus by the structure of a basic cycle we have $he_1 =
e_1h$. While $e_1$ is a generator of $G(C_1)$. Thus, $G(C_1)$ is
contained in the center of $G(C)$.

\vskip5pt

(ii) \ \ One can see this assertion from Theorem 5.1(ii) by
identifying $C_1$ with $C_d(n)$, and the claimed $g$ is exactly
$e_1$ in $C_d(n)$.

\vskip5pt

(iii) \ \ By Theorem 3.2 in [M2], as an algebra, $C$ is generated
by $C_1$ and $G(C)$.  By the proof of Theorem 3.1(ii)  $C_1$ is
generated by $g = e_1$ and a non-trivial $1, e_1$-primitive
element $x$, satisfying the given relation, together with

$$xe_1 = qe_1x$$
with $q$ a primitive $d$-th root of unity.

\vskip5pt

(iv) \ \  For any $h\in G$, since both $x\cdot h$ and $h\cdot x$
are non-trivial $h, gh$-primitive elements in $C$ (note $gh=hg$),
it follows that there exists $K$-functions $\chi$ and $\chi'$ on
$G$ such that

$$x\cdot h =\chi(h) h\cdot x + \chi'(h) (1-g)h.$$

\vskip5pt

We claim that $\chi$ is a one-dimensional representation of $G$
and $\chi'=0$.

\vskip5pt

By $x\cdot (h_1\cdot h_2) = (x\cdot h_1)\cdot h_2$, one infers
that

$$\chi(h_1\cdot h_2) = \chi(h_1)\chi(h_2)$$
and

$$\chi'(h_1\cdot h_2)=\chi(h_1)\chi'(h_2)+ \chi'(h_1).$$

Since $\chi(g) = q$ and $\chi'(g)=0$, it follows that
$\chi'(h\cdot g)=\chi'(h)$ for all $h \in G$. Thus, we have

$$\chi'(h) = \chi'(h\cdot g) = \chi'(g\cdot h) = \chi(g)\chi'(h),$$
which implies $\chi'=0$.

\vskip5pt

Since $x^d = \mu(1-g^d)$, it follows that one can make a choice
such that $\mu = 0$ if $d = n$. By $x^d \cdot h= \chi^d(h) h\cdot
x^d$ and $x^d=\mu (g^d-1)$ we see $\chi^d=1$ if $\mu \neq 0$ and
$g^d \neq 1$. $\blacksquare$

\vskip10pt

In order to classify non-semisimple monomial Hopf $K$-algebras, we
introduce the notion of group data.

\vskip10pt

\begin{defn} \  A group datum  $\alpha = (G, g, \chi, \mu)$ over $K$ consists of

\vskip5pt

(i) \ \ a finite group $G$, with an  element $g$ in its center;

\vskip5pt

(ii) \ \ a one-dimensional $K$-representation $\chi$  of $G$; and

\vskip5pt

(iii) \ \ an element $\mu \in K$, such that $\mu = 0$ if $o(g) =
o(\chi(g))$, and that if $\mu\ne 0$ then $\chi^{o(\chi(g))} = 1$.
\end{defn}

\vskip10pt

\begin{defn} \ Two group data $\alpha=(G, g,
\chi, \mu)$ and $\alpha'=(G', g', \chi', \mu')$ are said to be
isomorphic, if there exists a group isomorphism $f:
G\longrightarrow G^{'}$ and some $0\ne\delta \in K$ such that
$f(g) = g'$, $\chi = \chi' f$ and $\mu = \delta^d \mu'$.
\end{defn}

\vskip10pt

Lemma 5.2 permits us to introduce the following notion.

\vskip10pt

\begin{defn}\ Let $C$ be a non-semisimple monomial Hopf algebra.
A group datum $\alpha = (G, g, \chi, \mu)$ is called an induced
group datum of $C$ provided that

\vskip5pt

(i) \ \ $G = G(C)$;

\vskip5pt

(ii) \ \ there exists a non-trivial $1, g$-primitive element $x$
in $C$ such that

$$x^d=\mu (1-g^d), \ \ \ xh= \chi(h)hx,
\ \ \forall h \in G,$$ where $d$ is the multiplicative order of
$\chi(g)$.
\end{defn}

\vskip10pt

For example, $(Z_n,  \bar 1, \chi, \mu)$ with $\chi(\bar 1) = q$ is an induced group datum of
$A(n, d, \mu, q)$ (as defined in 1.7).

\vskip10pt

\begin{lem} \ \ (i) \ \ Let $C_1, C_2$ be
non-semisimple monomial Hopf algebras, $f: \ C\longrightarrow C'$
a Hopf algebra isomorphism,  and $\alpha = (G, g, \chi, \mu)$ an
induced group datum of $C$. Then $f(\alpha) = (f(G), f(g), \chi
f^{-1}, \mu)$ is an induced group datum of $C'$.

\vskip5pt

(ii) \ \ If $\alpha = (G, g, \chi, \mu)$ and $\beta = (G', g',
\chi', \mu')$ both are induced group data of a non-semisimple
monomial Hopf algebra $C$, then $\alpha$ is isomorphic to $\beta$.

\vskip10pt

Thus, we have a map $\alpha$ from the set of the isoclasses of
non-semisimple monomial Hopf $K$-algebras to the set of the
isoclasses of group data over $K$, by assigning each
non-semisimple monomial Hopf algebra $C$ to its induced group
datum $\alpha (C)$.
\end{lem}

\vskip5pt

\noindent{\bf Proof} \ \  The assertion (i) is clear by
definition.

(ii) \ \ By definition we have $G = G(C) = G'$. By definition
there exists a non-trivial $1, g$-element $x$, and also a
non-trivial $1, g'$-element $x'$. But according to Theorem 5.1
(ii) such $g$ and $g'$ turn out to be unique, i.e. $g = g' = e_1$
if we identify $C_1$ with $C_d(n)$. And according to the coalgebra
structure of $C$, and of $C_1\cong C_d(n)$, we have

$$x = \delta x'+ \kappa (1-g)$$
for some $\delta \neq 0, \kappa \in K$. It follows that

$$x\cdot h = \chi(h)h\cdot x = \chi(h) \delta h\cdot x' +
\chi(h)\kappa h\cdot (1-g)$$ and
$$x\cdot h = (\delta x' + \kappa (1-g))\cdot h =
\delta \chi'(h) h\cdot x' + \kappa h\cdot (1-g)$$
and hence $\chi=\chi'$ and $\kappa =0$. Thus

$$\mu(1-g^d) = x^d = (\delta x')^d = \delta^d \mu'(1-g^d),$$
i.e., $\mu = \delta^d\mu'$, which implies that $\alpha$ and
$\beta$ are isomorphic. $\s$

\vskip10pt

\sub{}\rm  For a group datum $\alpha = (G, g, \chi, \mu)$ over
$K$, define $A(\alpha)$ to be an associative algebra with
generators $x$ and all $h \in G $, with relations

\vskip10pt

\begin{center}
$x^d= \mu (1-g^d), \ \ \ xh = \chi(h) hx, \ \ \forall h\in G,$
\end{center}
where $d = o(\chi(g))$. One can check that ${\rm{dim}}_K A(\alpha)
= |G|d$ \ by Bergman's diamond lemma in [B] (here the condition "$\chi^d=1$
if $\mu \neq 0$" is needed). Endow $A(\alpha)$ with
comultiplication $\Delta$, counit $\epsilon$, and antipode $S$ by

\vskip10pt

\begin{center} $\Delta (x)=g \otimes x+ x\otimes 1 ,\ \
\epsilon(x)=0,$
\end{center}

\vskip10pt

\begin{center}
$\Delta (h)= h \otimes h,\ \ \epsilon(h)=1  ,   \ \ \forall h \in
G,$
\end{center}

\vskip10pt

\begin{center}
$ \ \ S(x)=g^{-1}x,  \ \ S(h)=h^{-1}, \ \, \forall h \in G.$
\end{center}

\vskip10pt

\noindent  It is straightforward  to verify that $A(\alpha)$ is
indeed a Hopf algebra.

\vskip10pt

\begin{lem} \ (i) \ \ For each group datum $\alpha = (G, g, \chi, \mu)$ over $K$,
$A(\alpha)$ is a non-semisimple monomial Hopf $K$-algebra, with
the induced group datum $\alpha$.

\vskip5pt

(ii) \ \ If $\alpha$ and $\beta$ are isomorphic group data, then
$A(\alpha)$ and $A(\beta)$ are isomorphic as Hopf algebras.

\vskip10pt

Thus, we have a map $A$ from the set of the isoclasses of group
data over $K$ to the set of the isoclasses of non-semisimple
monomial Hopf $K$-algebras, by assigning each group datum $\alpha$
to $A(\alpha)$.
\end{lem}

\vskip5pt

\noindent{\bf Proof} \ \ (i) \ \ Since ${\rm{dim}}_k A(\alpha)=
|G|d$, it follows that $\{hx^i \ |\ h\in G, \ i\leq d\}$ is a
basis for $A(\alpha)$. Let $\{a_1=1, \cdots, a_l\}$ be a set of
representatives of cosets
 of $G$ respect to $G_1$. For each $1 \leq i \leq l$,
let $A_i$ be the $K$-span of the set $\{a_ig^jx^k\ | \ 0 \leq j
\leq {n-1},\ 0\leq k \leq {d-1}\}$, where $n = |G_1|$. It is
straightforward to verify that

$$A(\alpha) = A_1\oplus \cdots \oplus A_l$$
as a coalgebra, and $A_i \cong A_j$ as coalgebras for all $1 \leq
i, \ j \leq l$. Note that there is a coalgebra isomorphism
$A_1\cong C_d(n)$, by sending $g^ix^j$ to $(j!_q)p_i^j$, where
$p_i^j$ is the path starting at $e_i$ and of length $j$.  This
proves that

 $$A(\alpha)\cong C_d(n)\oplus \cdots \oplus C_d(n)$$
as coalgebras.

\vskip5pt

(ii) \ \ Let $\alpha = (G, g, \chi f, \delta^d \mu)\cong \beta =
(f(G), f(g), \chi, \mu)$ with a group isomorphism $f:\
G\longrightarrow G'$. Then $F: \ A(\alpha)\longrightarrow
A(\beta)$ given by $F(x) = \delta x', \ \ F(h) = f(h), \ \ h\in
G$, is a surjective algebra map, and hence an isomorphism by
comparing the $K$-dimensions. This is also a coalgebra map, and
hence a Hopf algebra isomorphism. $\s$

\vskip10pt

The following theorem gives a classification of non-semisimple,
monomial Hopf $K$-algebras via group data over $K$.

\vskip10pt

\begin{thm} \ \ The maps $\alpha$ and $A$ above gives a one to one
correspondence between  sets

$$\{ \ \mbox{the isoclasses of non-semisimple monomial Hopf}
\ K\mbox{-algebras}\ \}$$ and

$$\{ \ \mbox{the isoclasses of group data over}\ K\ \}.$$
\end{thm}

\vskip5pt

\noindent {\bf Proof}\ \ By Lemmas 5.6 and 5.8, it remains to
prove that $C\cong A(\alpha(C))$ as Hopf algebras, which are
straightforward by constructions. $\blacksquare$

\vskip 10pt

\sub{}\rm\  A group datum $\alpha =(G, g, \chi, \mu)$ is said to
be trivial, if $G = \langle g\rangle \times N$, and the
restriction of $\chi$ to $N$ is trivial.

\vskip10pt

\noindent {\bf Corollary.} \ \ Let $\alpha =(G, g, \chi, \mu)$ be
a group datum over $K$. Then $A(\alpha)$ is isomorphic to $A(o(g),
o(\chi(g)), \mu, \chi(g)) \otimes KN$ as Hopf algebras, if and
only if $\alpha$ is trivial with $G = \langle g\rangle \times N$,
where $A(o(g), o(\chi(g)), \mu, \chi(g))$ is as defined in 1.7.

\vskip5pt

\noindent {\bf Proof} \ \  If $\alpha$ is trivial with $G =
\langle g\rangle \times N$, then
$$\alpha(A(o(g), o(\chi(g)), \mu, \chi(g))\otimes KN) = \alpha,$$
it follows from Theorem 5.9 that $$A(\alpha) \cong A(o(g),
o(\chi(g)), \mu, \chi(g))\otimes KN.$$

Conversely, we then have
$$\alpha = \alpha(A(\alpha)) = \alpha(A(o(g),
o(\chi(g)), \mu, \chi(g))\otimes KN)$$ is trivial. $\blacksquare$

\vskip10pt

\begin{rem} It is easy  to determine the automorphism group of $A(\alpha)$ with $\alpha=(G, g, \chi, \mu)$:
it is $K^* \times \Gamma $ if $\mu=0$, and $Z_d \times \Gamma$ if
$\mu \neq 0$, where $\Gamma:= \{f \in Aut(G) \ | \ f(g)= g, \ \chi
f = \chi\ \}$. \end{rem}

\vskip10pt

{\bf Acknowledgement:} We thank the referee for the helpful
suggestions.

\bibliography{}

\begin{thebibliography}{99}


\bibitem[ARS]{ARS} M. Auslander, I. Reiten, and S.O. Smal$\phi$,
Representation Theory of Artin Algebras, Cambridge Studies in Adv.
Math. 36, Cambridge Univ. Press,  1995.

\bibitem[AS]{AS}
N. Andruskiewitsch, and H. J. Schneider, Lifting of quantum linear
spaces and pointed Hopf algebras of order $p^3$.  J. Algebra
209 (1998), 658-691.

\bibitem[B]{B} G. M. Bergman, G.M., The diamond lemma for ring theory,
Advances in Mathematics 29 (1978) 178-218.

\bibitem[CM]{CM} W. Chin, and S. Montgomery, Basic coalgebras, Modular interfaces (Riverside, CA, 1995), 41-47, AMS/IP
Stud. Adv. Math. 4, Amer. Math. Soc., Providence, RI, 1997.

\bibitem[C]{C} C. Cibils, A quiver quantum group, Comm. Math.
Phys. 157 (1993), 459-477.

\bibitem[CR1]{cr4} C. Cibils, and M. Rosso, Algebres des chemins
quantique. Adv. Math. 125 (1997), 171-199 .

\bibitem[CR2]{CR} C. Cibils, C. and M. Rosso, Hopf quivers,
J. Algebra  254 (2) (2002), 241-251.

\bibitem[DK]{DK}  Yu. A. Drozd, and V. V. Kirichenko, Finite Dimensional
Algebras, Springer-Verlag, Berlin, Heidelberg, New York, Tokyo
1993.

\bibitem[GS]{ES} E. L. Green, and $\O$. Solberg, Basic Hopf algebras and
quantum groups, Math. Z. 229 (1998), 45-76.

\bibitem[HCZ]{HCZ} H. L. Huang, H. X. Chen, and P. Zhang, Generalized  Taft
Algebras. to appear in Algebra Colloquium.

\bibitem[K]{Kas} C. Kassel, Quantum Group. Graduate Texts in
Math. 155, Springer-Verlag, New York, 1995.

\bibitem[M1]{M} S. Montgomery, Hopf Algebras and Their Actions
on Rings, CBMS Regional Conf. Series in Math. 82, Amer. Math.
Soc., Providence, RI, 1993.

\bibitem[M2]{M2} S. Montgomery, Indecomposable coalgebras, simple comodules and
pointed Hopf algebras, Proc. of the Amer. Math. Soc. 123 (1995),
2343-2351.

\bibitem[Rad]{Rad} D. E. Radford, On the coradical of a finite-dimensional Hopf algebra, Proc. Amer. Math. Soc. 53(1) (1975), 9-15.

\bibitem[Rin]{Rin} C. M. Ringel, Tame Algebras and Integral Quadratic Forms, Lecture Notes in Math. 1099,
Springer-Verlag, 1984.

\bibitem[S]{Sw} M. E. Sweedler, Hopf Algebras. New York: Benjamin, 1969.

\bibitem[T]{Ta} E. J. Taft, The order of the antipode of finite
dimensional Hopf algebras, Proc. Nat. Acad. Sci. USA  68 (1971),
2631-2633.

\end{thebibliography}

\end{document}